\documentclass[12pt]{article}
\usepackage{geometry}
\usepackage{subfigure}
\usepackage{physics}
\usepackage{multirow}
\usepackage{booktabs}
\makeatletter
\newcommand\ackname{Acknowledgements}
\if@titlepage
\newenvironment{acknowledgements}{
	\titlepage
	\null\vfil
	\@beginparpenalty\@lowpenalty
	\begin{center}
		\bfseries \ackname
		\@endparpenalty\@M1
\end{center}}
{\par\vfil\null\endtitlepage}
\else

\usepackage{mathrsfs}
\usepackage{eucal}
\usepackage{amssymb,amsmath}
\usepackage{pgf}
\usepackage{tikz}
\usepackage{cite}
\usepackage{graphicx}
\usepackage{color}
\usepackage{amssymb}
\usepgflibrary{arrows}
\usetikzlibrary{calc,angles,positioning,intersections,quotes,backgrounds,decorations.pathreplacing}
\usetikzlibrary{decorations.markings}
\usepackage{amsmath}
\usepackage{amsthm}
\usepackage{accents}

\theoremstyle{remark}
\theoremstyle{definition}

\usepackage{comment}
\begin{document}
	\title{\textbf{One-parameter generalised Fisher information matrix: One random variable}}
	\author{Worachet Bukaew$^{\dagger,*,1}$ and  Sikarin Yoo-Kong$^{\dagger,2}$ \\	
		\small  \emph{$^\dagger$The Institute for Fundamental Study (IF),} \\
		\small\emph{Naresuan University, Phitsanulok, Thailand, 65000.}\\
		\small\emph{$^*$Center of Excellence in Theoretical and Computational Science (TaCs-CoE),}\\
		\small\emph{Faculty of Science, King Mongkut's University of Thonburi (KMUTT),}\\
		\small\emph{126 Pracha Uthit Rd., Bang Mod, Thung Khru, Bangkok, Thailand, 10140.} \\
		\small\emph{$^1$e-mail:warachetb62@nu.ac.th, $^2$e-mail:sikariny@nu.ac.th} 		
	}
	\date{}
	\maketitle
	\abstract
	We propose a generalised Fisher information or a one-parameter extended class of the Fisher information for the case of one random variable. This new form of the Fisher information is obtained from the intriguing connection between the standard Fisher information and the variational principle together with the non-uniqueness property of the Lagrangian. A generalised Cramér-Rao  inequality is also derived and a Fisher information hierarchy is also obtained from the two-parameter Kullback-Leibler divergence. An interesting point is that the whole Fisher information hierarchy, except for the standard Fisher information, does not follow the additive rule. Furthermore, the idea can be directly extended to obtain the one-parameter generalised Fisher information matrix for the case of one random variable with multi-estimated parameters. The hierarchy of the Fisher information matrices is obtained. The geometrical meaning of the first two matrices in the hierarchy is studied through the normal distribution. What we find is that these first two Fisher matrices give different nature of curvature on the same statistical manifold for the normal distribution.
	\\
	\\
	\textbf{Keywords}: Fisher information, Fisher information matrix, $q$-deformation, Tsallis index.
	
	\section{Introduction}\label{section1}
	There is no doubt that we are now in the ``information era". The information is physical \cite{InPh} and plays an essential role in modern physics ranging from thermodynamics, statistical mechanics, quantum mechanics to relativity. The birth of information theory can be traced back to the seminal paper of Shannon \cite{SHA2} on communication. The key quantity in this context is the entropy or more precisely ``Shannon entropy" as the mean value of information or uncertainty inherent in the possible outcomes. An interesting point is that the Shannon entropy takes the same form as the Gibbs-Boltzmann entropy in the context of statistical mechanics, which measures the configuration of the microscopic states, if one ignores the Boltzmann's constant. However, the notion of entropy was first introduced in the context of thermodynamics through the second law of thermodynamics-which is a bit more abstract relating to the heat flow in or out and the temperature of the system. However, if we trace back long before the breakthrough work of Shannon, Fisher purposed another information quantity, later known as Fisher information \cite{OrFis}, as an uncertainty measurement on estimating unknown parameters in the system. This means that the Fisher information allows us to probe into the internal structure of the system. Therefore, Shannon entropy and Fisher information provides a complete description of the system in the sense that the Fisher information can give an insight of what the system is made of and the Shannon entropy gives the system behavior in the big picture. The thing is that the Shannon differential entropy and Fisher information are connected which was first observed by Kullback \cite{KL}. With the Kullback insight, the Fisher information matrix can be obtained from the second derivative of the Kullback-Leibler divergence(or the relative entropy).
	\\
	\\
	The Fisher information has a wide range of applications as follows. There was shown that the Fisher information satisfies the H-theorem for the diffusion equation and, under certain conditions, the Fisher information gives the relative information measure of  Kullback \cite{Arow1}. Moreover, it had been illustrated that the Fisher information can be used to consider the relation between the arrow of time and symmetries of the Fokker-Planck equation \cite{FKP}. The Fisher information also appears in context of thermodynamic, since the thermodynamics’ first law is obtained for Fisher’s information measure known as the Fisher’s thermodynamics first-law \cite{Fisa1L}. Furthermore, the connection between Schroedinger equation and Boltzmann transport equation can be established through the Fisher inforamtion \cite{FisaShr}. Another interesting phenomena is that, for a system obeying conservation of flow, the rate of increase of entropy partially depends on the Fisher information  \cite{FisaIne1}, see also \cite{FisaIne2}. This feature was coined as an addendum to the second law. 
	\\
	\\
	A generalised version of the Shannon entropy was first introduced by Renyi \cite{Renyi}. The Renyi entropy comes with a parameter and, with a suitable limit, the Shannon entropy can be recovered. Then one could think that the Renyi entropy is a one-parameter extended class of the Shannon entropy. On the statistical mechanical side, a generalised version of the Gibbs-Boltzmann entropy was proposed by Tsallis \cite{Tsallis}. Again, this Tsallis entropy comes with a parameter and the Gibbs-Boltzmann entropy can be recovered with a suitable limit. One main feature of the Tsallis entropy is the non-additive property directly related to the non-extensivity of the system. Consequently, this leads to a new kind of research area known as the Tsallis statistics with a wide range of applications in statistics, physics, chemistry, economics, and biosciences \cite{BTsallis,TsBook}. Recently, an application in astrophysic for Tsallis entropy is studied \cite{Almeida}, since the Boltzmann–Gibbs entropy is no longer a good mathematical tool to describe the systems. On the other hand, several extensions of the Fisher information have been proposed with different aspects \cite{FracEn,GenFis,FracEn2,FracEn3,FANON1,FANON2} to serve different uses in statistics. Recently, the one-parameter generalised Fisher information in the case of one random variable and one estimated variable is derived \cite{WS}. Then in this contribution, we will extend the idea to the case of Fisher information matrix with multi-estimated parameters and one random variable. The key motivation and derivation are still based on the intriguing connection between Fisher information and variational principle observed by Frieden \cite{FD1,FD2,FD3} together with one-parameter extended class of the Lagrangian \cite{MulLag}.
	\\
	\\
	The remaining body of this work will be divided into two parts. The first half will be the overview of the work based on \cite{WS} with some extended points as follows. In section \ref{section2}, a brief review of the standard Fisher information together with notations will be given. In section \ref{section3}, one parameter extended class of the Fisher information is given by employing the connection with the variational principle. The Fisher information hierarchy will be therefore obtained. In section \ref{section4}, the extended Cramér-Rao inequality and non-additive property are given. In section \ref{section5}, the connection between two-parameter Kullback-Leibler divergence and Fisher information hierarchy will be established. In section \ref{section7}, which is a second half of this work, one parameter extended class of the Fisher information matrix will be obtained by directly employing the idea from the first half of the paper. The normal distribution will be used as the example to examine the geometrical nature of the Fisher information matrices. The last section will be about the conclusion and discussion. 
	
	\section{Fisher information and Fisher information matrix}\label{section2}
	We shall start in this section to give a brief review of the basic Fisher information and the connection with Kullback-Leibler divergence. Given a set of random variables $X = (x_1, x_2, ..., x_N)$, which belongs in the sample space $\Omega$ $(X \in \Omega)$, and associated probability density functions $p(x_i\mid \theta)$, which is identical and independent, we define the likelihood as 
	$L(\theta\mid X)=\prod_{i=1}^N p(x_i\mid \theta)$, 
	where $\theta$ is an unknown parameter in the probability model. In principle, the likelihood measures how good the statistical model is comparing to the sample of data $X$ for given the values of the unknown parameter $\theta$. Then what we need is to find maximum value of the likelihood $\partial  L(\theta\mid X) /	 \partial \theta  =0$.
	However, in practice, the likelihood might not be well behave in calculation. Since logarithm is monotonic increasing function, we shall work with the log-likelihood $\partial \log L(\theta\mid X) / \partial \theta  =0$ (also called the score function), which gives the same maximum points.
	What we are interested is the nature of the maximum peak as an indicator for how good it is for the estimated value. If the peak is thin, the estimated value is fairly determined. On the other hand, if the peak is broad, the estimated value is uncertainly determined. With this reason, we therefore consider the second derivative of the log-likelihood averaged over all possible random variables as the measure called the Fisher information
	\begin{eqnarray}
		I(\theta)\equiv - \left< \frac{\partial^2 }{\partial \theta^2}\log L(\theta \mid X)\right>\label{Fis1} \;.
	\end{eqnarray}
	Employing the fact that 
	\begin{eqnarray}
		\frac {\partial }{\partial \theta} \log L(\theta\mid X)=\frac{1}{L(\theta\mid X)}\frac{\partial L(\theta\mid X)}{\partial \theta}\;,
	\end{eqnarray}	
	and
	\begin{eqnarray}	
		\left< \frac{\partial }{\partial \theta}\log L(\theta \mid X)\right>=0\;.
	\end{eqnarray}
	The Fisher information can be re-expressed as 
	\begin{eqnarray}
		I(\theta)= - \left< \frac{\partial^2 }{\partial \theta^2}\log L(\theta\mid X)\right> 
		=\left< \left( \frac{\partial  }{\partial \theta}\log L(\theta\mid X) \right)^2 \right>  
		= \text{Var}\left[ \frac{\partial  }{\partial \theta}\log L(\theta\mid X) \right] \;.
	\end{eqnarray}
	Here, we can conclude that Fisher information tells us how much we know about the internal structure from data with a given space of outcomes (sample space) $\Omega$. 	
	\\
	\\
	The Fisher information also provides a information lower bound on the variance of an unbiased estimator for a parameter. This relation is known as the Cramér-Rao inequality. To obtain such relation, one can start to consider the unbiased estimator given by
	\begin{eqnarray}
		B(\hat{\Theta})\equiv\left<\hat{\Theta}-\theta \right>=\int_{\Omega}...\int_{\Omega} (\hat{\Theta}-\theta) L(\theta\mid X) dX=0\;,\label{Beq}
	\end{eqnarray}
	where $\hat{\Theta}=h(x_1, x_2, ..., x_N)$ is a point estimator. From now on, we might neglect subscription $\Omega$ for the sake of simplicity. Next, we consider the derivative of \eqref{Beq} with respect to the parameter $\theta$ resulting in
	
	\begin{eqnarray}
		\frac{\partial}{\partial \theta} \left<\hat{\Theta}-\theta \right>&=& - \int...\int L(\theta\mid X) dX+ \int...\int (\hat{\Theta}-\theta) \frac{\partial}{\partial \theta} L(\theta\mid X) dX \;. 
	\end{eqnarray}
	Using the fact that $ \int...\int L(X \mid \theta ) dX=1 $, we obtain 
	\begin{eqnarray}
		\int...\int \left[ (\hat{\Theta}-\theta) \cdot L^{1/2}(\theta\mid X) \right] \left[ \left(\frac{\partial}{\partial \theta} \log L(\theta\mid X) \right) \cdot L^{1/2}(\theta\mid X) \right] dX &=& 1 \; .
	\end{eqnarray}
	Applying Cauchy–Schwarz inequality, we obtain
	\begin{eqnarray}
		\frac{1}{\int...\int  {(\hat{\Theta}-\theta)^2  L(\theta\mid X)}dX}   &\leq&  \int...\int \left(\frac{\partial}{\partial \theta} \log L(\mathbf{x} \mid \theta )\right)^2  L(\theta\mid X) dX \nonumber\\
		\frac{1}{ \text{Var}(\hat{\Theta})}  &\leq&  I(\theta)\;.\label{CRB}
	\end{eqnarray}
	What we have in \eqref{CRB} is that the variance of any such estimator is at least as much as the inverse of the Fisher information.
	\\
	\\
	In general, the estimated parameters could come in a set i.e., $\theta=(\theta_1,\theta_2,...,\theta_n)$. What we have now is the element $(i,j)$ of Fisher information matrix: $[I(\theta)]_{ij}$ given by
	\begin{eqnarray}
		[I(\theta)]_{ij}= -\left< \frac{\partial^2}{\partial\theta_i\partial\theta_j} \log L(\theta|X)\right>=\left< \left( \frac{\partial  }{\partial \theta_i}\log L(\theta\mid X) \right) \left( \frac{\partial  }{\partial \theta_j}\log L(\theta\mid X) \right) \right>\;. \label{Fismatrix}
	\end{eqnarray}
	The interesting point is that the Fisher matrix can also be obtained by considering the relative entropy or Kullback-Leibler divergence between two distribution $p(X)$ and $q(X)$ on the probability manifold
	\begin{equation}
		KL(p||q)=\int...\int p(X)\log\left(\frac{p(X)}{q(X)}\right)dX\;.
	\end{equation}
	Then the Kullback-Leibler divergence between two probability distributions  $L(\theta|X)$ and $L(\theta'|X)$, parametised by $\theta$, is given by
	\begin{equation}
		D(\theta,\theta')\equiv KL(L(\theta|X)||L(\theta'|X)) =\int...\int L(\theta|X)\log\left(\frac{L(\theta|X)}{L(\theta'|X)}\right)dX\;.
	\end{equation}
	For $\theta$ being fixed, the Kullback-Leibler divergence can be expanded around $\theta$ as
	\begin{equation}
		D(\theta,\theta')=\frac{1}{2}\left(\theta'-\theta\right)^{\text T}\left(\frac{\partial^2}{\partial \theta'_i\partial \theta'_j}D(\theta,\theta') \right)\Big\vert_{\theta=\theta'}\left(\theta'-\theta\right)+\mathcal{O}\left(\left(\theta'-\theta\right)^2\right)\;,
	\end{equation}
	where the second order derivative is 
	\begin{equation}
		\left(\frac{\partial^2}{\partial \theta'_i\partial \theta'_j}D(\theta,\theta') \right)\Big\vert_{\theta=\theta'}=-\int...\int\left( \frac{\partial^2}{\partial \theta'_i\partial \theta'_j}L(\theta'|X)\right)\Big\vert_{\theta'=\theta}L(\theta|X)dX=[I(\theta)]_{ij}\;.
	\end{equation}
	With this connection, one may intuitively interpret the Fisher information as the metric between two point on the probability manifold. However, the Kullback-Leibler divergenece is not symmetric and does not follow the triangle inequality \cite{Stone}. Then the Fisher information cannot be treated as a true metric.
	
	\section{One-parameter generalised Fisher information: One random variable}\label{section3}
	In this section, we shall derive a new type of the Fisher information called a one parameter generalised Fisher information. One of the key features which will play an important role to obtain such new Fisher information is the variational principle. The intriguing connection between the variational principle and Fisher information was proposed by Frieden \cite{FD3}. Another feature is the non-uniqueness property of the Lagrangian which leads to the non-standard form of the Lagrangian called the multiplicative form \cite{MulLag}.  
	\\
	\\
	Let us now define $\theta$ as actual value of a measurement quantity, $X=(x_1,x_2,..,x_N)$ with $N$ outcomes and $Y=(y_1,y_2,...,y_N)$ as random errors associated with each measurement. Then we have 
	
	\begin{equation}
		x_i = \theta + y_i \;. \label{mesur1}
	\end{equation}
	Next, let $p(x_i|\theta)$ be a probability distribution and $L(\theta\mid X)=\prod_{i=1}^N p(x_i\mid \theta)$ be likelihood of probability distribution over the $x$'s with respect to $\theta$. Recalling the Fisher information \footnote{Here, we prefer the natural logarithm function.} \eqref{Fis1}, we can write the Fisher information in additive (see more details in \cite{FD3}) form as 
	\begin{align}
		I(\theta)= \int \left( \frac{\partial }{\partial \theta}\ln L(\theta\mid X)\right)^2 L(\theta\mid X)\prod_{i=1}^N dx_i
		=\sum_{i=1}^{N} \int \left( \frac{\partial }{\partial \theta}\ln p(x_i\mid \theta)\right)^2 p(x_i\mid \theta) dx_i\;.
	\end{align}
	In the case $N=1$, we do have $p(x \mid \theta )=p(x-\theta )=p(y)$. The  Fisher information is simply reduced to
	\begin{equation}
		I[p(y)] = \int \frac{\left( p'(y)\right)^2}{p(y)}dy\;\;\; \text{where}\;\;\; p'(y)=\frac{dp(y)}{d y}. \label{fisher2}
	\end{equation} 
	Next, we do a transformation such that $q(y)=\sqrt{p(y)}$ resulting in
	\begin{equation}
		I[q(y)] = 4 \int q'^2(y) dy, \;\; \text{where} \;\;q'(y)= \frac{dq(y)}{d y}.\label{F1}
	\end{equation}
	We find that the Fisher information is now a functional with the input function $q(y)$. Here comes to an interesting point. If we define $I[q]\equiv S[q]$ as an action functional and $\mathscr{L}(q',q;y)\equiv 4 q'^2(y)$ as the Lagrangian, the variational principle would give the Euler-Lagrange equation
	
	\begin{equation}
		\frac{\partial \mathscr{L}(q',q;y)}{\partial q(y)}- \frac{d}{dy} \left( \frac{\partial \mathscr{L}(q',q;y)}{\partial q'(y)}\right)=0,
	\end{equation}
	resulting in $-8q''(y)=0$. This second order differential equation describes the motion for a free particle (how a position $q$ change with time $y$). This would mean that the Fisher information \eqref{F1} could be remarkably treated as the action functional for the free particle and, of course, in the absence of the interaction, the equation of motion in physics is a direct result of extramising the Fisher information: $\delta I[q]=0$. 
	\\
	\\
	Next, what we would like to focus is the property of the Lagrangian known as the non-uniqueness property. Commonly, two Lagrangians, differing by the total derivative with respect to time of some function $F(q,y)$, would give an identical equation of motion on extreamising the action. However, one could ask an inverse question as follows. Imposing the equation of motion, could we solve all possible Lagrangians directly from the Euler-Lagrange equation? The answer is definitely yes and this problem is known as the inverse problem of the calculus of variations\cite{inverv}. Recently, Sarawuttinack et al \cite{MulLag} proposed a new type of Lagrangian called the multiplicative form for the case of one degree of freedom. Here we shall employ the same technique and propose an alternative Lagrangian for $\mathscr{L}(q',q;y)=4q'^2(y)$ as $\mathscr{L}_\lambda(q',q;y)=4[e^{\lambda q'^2(y)}-1]/\lambda$, where the parameter $\lambda$ has to be an inverse square of the velocity. Of course, the Lagrangian $\mathscr{L}_\lambda(q',q;y)$ can be treated as one-parameter extended class of the standard Lagrangian $\mathscr{L}(q',q;y)$.  It is not difficult to see that these two Lagrangians give exactly the same equation of motion. By considering the limit $\lambda\rightarrow 0$, one finds that $\lim\limits_{\lambda \to 0}\mathscr{L}_\lambda(q',q;y)=\mathscr{L}(q',q;y)=4q'^2(y)$. Then what we have is the action functional in the form
	\begin{equation}
		I_\lambda[q(y)] = \frac{4}{\lambda} \int \bigg[e^{\lambda q'^2(y)}-1\bigg] dy\;.\label{LF} 
	\end{equation}
	We shall call \eqref{LF} as a one parameter generalised Fisher information. The reason can be seen as follows. If we expand the functional \eqref{LF} with respect to the parameter $\lambda$, we obtain 
	
	\begin{eqnarray}
		I_\lambda[q(y)]  &=& 4 \int q'^2(y) dy + 4\frac{\lambda}{2!} \int q'^4(y)dy + 4\frac{\lambda^2}{3!} \int q'^6(y)dy+...\nonumber  \\
		&=& I_1[q(y)]+  \frac{\lambda}{2!} I_2[q(y)] + 	\frac{\lambda^2}{3!} I_3[q(y)] +... \;.\label{F2}
	\end{eqnarray}
	What we have in \eqref{F2} is a hierarchy $\{I_1,I_2,I_3,...\}$, where the first three functionals are 
	\begin{subequations}\label{FFH}
		\begin{eqnarray}
			I_1 &=& 4  \int\left(\frac{p'(y)}{2q(y)}\right)^2 dy= \frac{4}{{2^2}}  \int \frac{p'^2(y)}{p(y)} dy = \frac{4}{{2^2}} \int \left( \frac{\partial}{\partial \theta} \ln p(x|\theta) \right)^2 p(x|\theta) dx\;,\label{FFF11}\\
			I_2 &=& 4  \int\left(\frac{p'(y)}{2q(y)}\right)^4 dy= \frac{4}{{2^4}}  \int \frac{p'^4(y)}{p^2(y)} dy = \frac{4}{{2^4}} \int \left( \frac{\partial}{\partial \theta} \ln p(x|\theta) \right)^4 p^2(x|\theta) dx\;,\label{FFF1}\\
			%
			I_3& =& 4  \int\left(\frac{p'(y)}{2q(y)}\right)^6 dy= \frac{4}{{2^6}}  \int \frac{p'^6(y)}{p^3(y)} dy = \frac{4}{{2^6}} \int \left( \frac{\partial}{\partial \theta} \ln p(x|\theta) \right)^6 p^3(x|\theta) dx\;.
		\end{eqnarray} 
	\end{subequations}
	The first term is nothing but the standard Fisher information $I_1[q]=I[q]$ coinciding with the limit $\lim\limits_{\lambda \to 0} I_\lambda[q(y)]  = I[q(y)]$. With the structure in \eqref{FFH}, one could easily write $I_n [p(y)]$ as
	\begin{eqnarray}
		I_n (\theta) = \frac{4}{{2^{2n}}} \int \left( \frac{\partial}{\partial \theta} \ln p(x|\theta) \right)^{2n} p^n(x|\theta) dx\;,\;\; n=1,2,3,...\;\;\;. \label{extenF}
	\end{eqnarray} 
	Next point is that the standard Fisher information is in the average of the score function but the rest in the hierarchy is not. Then we shall seek a transformation to express the higher order Fisher information in the statistical average. We shall first consider the second function in \eqref{FFF1} and introduce a new variable $\phi_1 \equiv p^2$ such that $p'(y) = \phi'_1(y)/2p(y)$, 
	resulting in	
	\begin{eqnarray}
		I_2 (\theta) =  \frac{4}{4^4} \int \frac{ \phi'^4_1(y)}{\phi^4_1(y)} \phi_1(y) dy,\;\;\;
		\text{or}
		\;\;\;
		I_2(\theta) = \frac{4}{4^4}  \int \left[ \frac{\partial}{\partial \theta}\ln \phi_1(x|\theta) \right]^4 \phi_1(x|\theta)dx\;.\label{LF2}
	\end{eqnarray}	
	We shall call \eqref{LF2} as the 2$^\text{nd}$ order Fisher information. 
	We can proceed the same technique of transformation, then the $n^{\text{th}}$ order Fisher information \eqref{extenF} becomes
	\begin{eqnarray}
		I_n(\theta )&=& \frac{4}{(2n)^{2n}} \int \left[ \frac{\partial}{\partial \theta}\ln \phi_{n-1}(x|\theta) \right]^{2n} \phi_{n-1}(x|\theta)dx, \label{GnEF}
	\end{eqnarray}	
	where $\phi_{n-1}(y)=p^{n}(y)$ and the generalised Fisher information \eqref{F2} can be expressed in terms of infinite series as 
	\begin{equation}
		I_\lambda(\theta )=\sum_{n=1}^\infty \frac{\lambda^{n-1}}{n!}I_n(\theta )\;.
	\end{equation}
	At this point, we may treat \eqref{F2} as the generating function for the entire hierarchy of the Fisher information by expanding with respect to the parameter $\lambda$. 
	\\
	\\
	\textbf{Working example}: The remaining thing is that what is the meaning of the higher order $n^{\text{th}}$ Fisher information? For simplicity, we shall try to interpret the 1$^\text{st}$ order and 2$^\text{nd}$ order of the Fisher information in the hierarchy. To make thing more transparent, we shall work with the well known example: the normal distribution with two parameter $\theta= \{ \mu,\sigma \} $
	\begin{eqnarray}
		p(x|\theta= \{ \mu,\sigma \}) = \frac{1}{\sqrt{2\pi \sigma^2} } e^{-\frac{(x-\mu)^2}{2\sigma^2}}\;,
	\end{eqnarray}	
	where $\mu$ is mean value and $\sigma^2$ is variance, respectively.
	Here, the parameter $\sigma$ is given and the parameter $\mu$ will be estimated. The standard Fisher information, the 1$^\text{st}$ term in the hierarchy, is $1/\sigma^2$. This means that for normal distribution, Fisher information is proportional to the variance of the density function. Therefore, when the normal distribution is widely range (our data is in sort of spread), that is high value of $\sigma$, then we have a low Fisher information. On the other hand, if normal distribution comes in a narrow range, that is low value of $\sigma$, then we have a high Fisher information. These features imply that Fisher information can be used to measure the degree of disorder of a system as same as the entropy can do \cite{FD1}.
	\\
	\\
	For the 2$^\text{nd}$ terms in the Fisher information heirarchy, we first make the transformation on the probability such that  $\phi_1 = p^2$, resulting in
	\begin{eqnarray}
		\phi_1(x|\theta= \{ \mu,\sigma \})=\frac{1}{2\pi \sigma^2} e^{-\frac{(x-\mu)^2}{\sigma^2}}, 
	\end{eqnarray}   	
	and
	\begin{eqnarray}
		\nonumber \ln \phi_1(x|\theta= \{ \mu,\sigma \})
		&=&\ln{\frac{1}{2 \pi \sigma^2}}-\frac{(x-\mu)^2}{\sigma^2}\\
		&=& \ln{\frac{1}{2 \pi \sigma^2}}-\frac{x^2-2\mu x+\mu^2}{\sigma^2}\;.
	\end{eqnarray}
	Then we can compute the log-likelihood by	
	\begin{eqnarray}
		\frac{\partial \ln \phi_1(x|\theta= \{ \mu,\sigma \})}{\partial \mu}= \frac{2(x-\mu)}{\sigma^2}\;.
	\end{eqnarray}	
	We then have
	\begin{align}
		I_{2}[\theta] =& \frac{4}{(4^4)}  \int_{-\infty}^{\infty} \left[\frac{\partial}{\partial \mu} \ln  \phi_1(x|\theta= \{ \mu,\sigma \})\right]^4 \phi_1(x|\theta= \{ \mu,\sigma \}) dx \nonumber\\
		=& \frac{2^4 4 }{(4^4)\sigma^8} \int_{-\infty}^{\infty} [x^4-4\mu x^3+6x^2 \mu^2-4 \mu^3 x+ \mu^4 ] \phi_1(x|\theta= \{ \mu,\sigma \}) dx \nonumber \\
		=& \frac{2^4 4 }{(4^4)\sigma^8}  \left[ m_4-4 \mu m_3+6 \mu^2 m_2-3 \mu^4 \right] \nonumber \\
		=& \frac{1}{4 (\sigma^2)^4} m_4'\;,
	\end{align}
	where $m_n = \left< x^n \right>$ is called the $n^{\text{th}}$ statistical raw moment and $m_4'=\left< (x-\mu)^4  \right>$ is the 4$^\text{th}$ central moment.
	We can see that the 2$^\text{nd}$ order Fisher information depends on the 4$^\text{th}$ central moment for probability $\phi_1(x|\mu,\sigma)$ and the reciprocal of the fourth power of the variance. The result seems to suggest that the 2$^\text{nd}$ order Fisher information tells us about the outliers and tailedness\cite{Probtail}. Here let us elaborate a bit more. Suppose there are two normal distributions with different variances. For the large variance, the distribution will be broad and of course this means that it has higher chance to get the relatively rare events or long tail. On the other hand, the distribution will be quite sharp for the small variance and the tail will be short. Then with this particular example, we can see that the 2$^\text{nd}$ order Fisher information can possibly used to indicate how good we could make the estimation of the unknown parameter through the nature of the tail of the distribution(which, of course, is related to the variance of the distribution)
	\section{Generalised Cramér-Rao inequatlity and non-additive property }\label{section4}
	Here in this section, we will construct the Cramér-Rao inequality associated with the Fisher information hierarchy given in the previous section. To achieve the goal, we start with 
	\begin{eqnarray}
		\left<\hat{\Theta}-\theta \right>=\int (\hat{\Theta}-\theta) p^q(x \mid \theta ) dx=0\;,
	\end{eqnarray}
	which is known as the $q$-expectation value \cite{Imcomp1,Imcomp2,Imcomp3} .
	Taking the $1^{\text{st}}$ derivative, we obtain
	\begin{eqnarray}
		\frac{\partial}{\partial \theta} \left<\hat{\Theta}-\theta \right>&=&\int \frac{\partial}{\partial \theta}(\hat{\Theta}-\theta) p^q(x \mid \theta ) dx+ \int (\hat{\Theta}-\theta) \frac{\partial}{\partial \theta} p^q(x \mid \theta ) dx \nonumber \\
		&=& -\int  p^q(x \mid \theta ) dx + q\int (\hat{\Theta}-\theta) p^{q-1}(x \mid \theta ) \frac{\partial}{\partial  \theta} p(x \mid \theta ) dx \nonumber \\ &=& -\int  p^q(x \mid \theta ) dx + q\int (\hat{\Theta}-\theta) p^{q-1}(x \mid \theta ) p(x \mid \theta ) \frac{\partial }{\partial  \theta}\ln p(x \mid \theta )   dx \nonumber \\&=& -Q_q+qJ=0 \label{QJ}\;,
	\end{eqnarray}
	where 
	\begin{eqnarray}
		Q_q&\equiv&\int  p^q(x \mid \theta ) dx\;,\\
		J&\equiv&\int (\hat{\Theta}-\theta) p^{q-1}(x \mid \theta ) p(x \mid \theta ) \frac{\partial  }{\partial  \theta} \ln p(x \mid \theta ) dx\;.
	\end{eqnarray}
	Conventionally, the term $Q_q$ is known as the information generating function \cite{Imcomp3} with Tsallis index $q$\cite{Tsallis}(we are now interested in case $q \geq 1$) and it is also called the incomplete normalization. The $p^q$ is called the effective probability\cite{Imcomp2}. Next, we rewrite the $J$ in the form
	\begin{eqnarray}
		J &=& \int \left[(\hat{\Theta}-\theta) \right] \left[ \bigg(\frac{\partial }{\partial  \theta} \ln p(x \mid \theta )\bigg)\; p^{q-1}(x \mid \theta )\right] p(x \mid \theta ) dx \;.\label{Je}
	\end{eqnarray}
	Applying the Hölder's inequality \cite{inequa} to \eqref{Je}, we obtain
	\begin{eqnarray}
		J &\leq& \left[\int (\hat{\Theta}-\theta)^{\beta} p(x \mid \theta ) dx \right]^{1/\beta} \left[\int \left(\frac{\partial  }{\partial  \theta}\ln p(x \mid \theta )\right)^{\alpha} (p^{q-1}(x \mid \theta ))^\alpha p(x \mid \theta ) dx \right]^{1/\alpha} \nonumber \\
		&=& \left[\int (\hat{\Theta}-\theta)^{\beta} p(x \mid \theta ) dx \right]^{1/\beta} \left[\int \left(\frac{\partial  }{\partial  \theta}\ln p(x \mid \theta ) \right)^{\alpha} p^{\alpha(q-1)+1}(x \mid \theta ) dx \right]^{1/\alpha} \label{IE1},
	\end{eqnarray}
	where Hölder conjugates $\alpha$ and $\beta $ are related with the condition $1/\alpha +1/\beta =1 $ for $\alpha, \beta = [1,\infty]$.
	Finally, employing \eqref{QJ}, the inequality \eqref{IE1} becomes
	\begin{eqnarray}
		\frac{Q_q}{q}=\frac{\int  p^q(x \mid \theta ) dx}{q} \leq \left[\int (\hat{\Theta}-\theta)^{\beta} p(x \mid \theta ) dx \right]^{1/\beta} \left[\int \left(\frac{\partial }{\partial  \theta}\ln p(x \mid \theta )\right)^{\alpha} p^{\alpha(q-1)+1}(x \mid \theta ) dx \right]^{1/\alpha}\;,\label{InqEX}
	\end{eqnarray}
	which is our generalised Carmer-Rao inequality. It is not difficult to see that if one takes $q=1$, $\beta=2$ and $\alpha=2$, the standard Carmer-Rao inequality can be recovered. For  $\alpha = 4$, $\beta=4/3$ and $q=5/4$, we obtain
	\begin{eqnarray}
		\frac{  4 \; Q_{5/4}}{5^4 \left< (\hat{\Theta}-\theta)^{4/3} \right>^{3}}&\leq& I_2\;,		
	\end{eqnarray}
	which is the Cramér-Rao inequality for the 2$^\text{nd}$ order Fisher information. Basically, the inequality \eqref{InqEX} provides the Cramér-Rao  bound for the whole Fisher information hierarchy as shown in table \ref{Ta1}.
	\begin{table}[htbp]
		\caption{The three parameters associated with the first four Carmer-Rao inequalities. }
		\begin{center}
			\begin{tabular}{|c|c|c|c|}
				\hline
				\multirow{2}{5cm}{\centering $n^\text{nd}$ Carmer-Rao inequality}& \multicolumn{3}{p{3cm}|}{\centering Parameters} \\
				\cline{2-4} & \multicolumn{1}{c|}{$q$} & \multicolumn{1}{c|}{$\beta$} & \multicolumn{1}{c|}{$\alpha$} \\ \hline
				$1^\text{st}$ order & 1 & 2 & 2 \\
				$2^\text{nd}$ order & 5/4 & 4/3 & 4 \\
				$3^\text{rd}$ order & 4/3 & 6/5 & 6 \\
				$4^\text{th}$ order & 11/8 & 8/7 & 8\\
				\hline
			\end{tabular}
		\end{center} \label{Ta1}
	\end{table} 
	\newpage
	
	\noindent	Next, we will investigate the additive property of the higher order Fisher information. For simplicity, we shall start with the 2$^\text{{nd}}$ order Fisher information. Suppose a system composed of two independent identically subsystems with random variables $X=(x_1, x_2)$, where superscription denote for subsystems. The joint probability of the two subsystems is given by $p_{12}\equiv p(x_1,x_2|\theta)=p(x_1|\theta) p(x_2|\theta)\equiv p_1p_2$. 
	What we have for the 2$^\text{{nd}}$ order Fisher information is
	\begin{align}
		I_2 (p_{12}) &= \frac{4 }{2^4 } \int \int \left( \frac{\partial}{\partial \theta}\ln (p_1 p_2) \right)^4 p^2_{1} p^2_{2}dx_1 dx_2 \nonumber \\
		&= \frac{4 }{2^4 }  \bigg[ \int \left(\frac{\partial}{\partial \theta}\ln p_{1}\right)^4p^2_{1} dx_1 \int p^2_{2}dx_2
		+ 4 \int \left(\frac{\partial}{\partial \theta}\ln  p_{1}\right)^3  p^2_{1} dx_1 \int \left(\frac{\partial}{\partial \theta} \ln p_{2}\right) p^2_{2}dx_2 \nonumber \\
		&+6\int \left(\frac{\partial}{\partial \theta}\ln  p_{1}\right)^2  p^2_{1} dx_1 \int \left(\frac{\partial}{\partial \theta} \ln p_{2}\right)^2 p^2_{2}dx_2 
		+ 4 \int \left(\frac{\partial}{\partial \theta}\ln  p_{1}\right)  p^2_{1} dx_1 \int \left(\frac{\partial}{\partial \theta} \ln p_{2}\right)^3 p^2_{2}dx_2 \nonumber \\
		&+ \int p^2_{1} dx_1 \int \left(\frac{\partial}{\partial \theta} \ln p_{2}\right)^4 p^2_{2}dx_2  \bigg] \nonumber \\ 
		&= \frac{4 }{2^4 }  \bigg[ Q_{2}(p_2)\int \left(\frac{\partial}{\partial \theta}\ln p_{1}\right)^4 p^2_{1} dx_1 
		+6\int \left(\frac{\partial}{\partial \theta}\ln p_{1}\right)^2  p^2_{1} dx_1 \int \left(\frac{\partial}{\partial \theta} \ln p_{2}\right)^2 p^2_{1}dx_2 \nonumber \\
		&+ Q_{2}(p_1) \int \left(\frac{\partial}{\partial \theta} \ln p_{2}\right)^4 p^2_{2}dx_2  \bigg] \nonumber \\
		&= \frac{1}{4}  \bigg[ Q_{2}(p_2) I_2(p_{1})+  Q_{2}(p_1) I_2(p_{2})+6 I(p_{1})I(p_{2})\bigg]\;.\label{2FF}
	\end{align}
	Here see that the 2$^\text{{nd}}$ order Fisher information does not follow the additive rule. With the result in \eqref{2FF}, it is not difficult now to see that the n$^\text{{th}}$ order Fisher information could give
	\begin{align}
		I_n (p_{12}) =& \frac{4 }{2^{2n} } \bigg[  {{2n}\choose{0}} Q_{n}(p(x_2|\theta))\int \left(\frac{\partial }{\partial \theta}\ln p(x_1|\theta) \right)^{2n} p^n(x_1|\theta) dx_1 \nonumber \\
		& + \sum_{k=2}^{2n-2} {{2n}\choose{k}} \int \left(\frac{\partial }{\partial \theta}\ln p(x_1|\theta) \right)^{2n-k} p^n(x_1|\theta) dx_1 \int \left(\frac{\partial }{\partial \theta} \ln p(x_2|\theta) \right)^{k} p^n(x_2|\theta) dx_2 \nonumber \\
		& + {{2n}\choose{2n}} Q_{n}(p(x_1|\theta)) \int \left(\frac{\partial }{\partial \theta}\ln p(x_2|\theta) \right)^{2n} p^n(x_2|\theta) dx_2 \bigg]\;,	
	\end{align}
	where the first and last terms refer to the Fisher information for each subsystem and the middle one is the crossing-term. Therefore, our Fisher information hierarchy does not follow the additive property, except for the standard Fisher information: $n=1$.
	%
	%
	%
	%
	%
	\section{The Kullback–Leibler divergence revisited}\label{section5}
	Here in this section, we shall investigate on the connection between our Fisher information hierarchy and the Kullback–Leibler divergence. We shall begin with the Kullback–Leibler divergence	
	\begin{eqnarray}
		\mathit{D}(p(y)\parallel q(y) ) =\int p(y)\ln\left(\frac{p(y)}{q(y)}\right)dy\;,
	\end{eqnarray}
	where $p$ and $q$ are two different points on the probability manifold. If $q(y)=p(y+\Delta)=p(y)+\Delta p'(y)$, where $p'(y)=dp(y)/dy$, we have
	\begin{eqnarray}
		\mathit{D}(p(y)\parallel p(y)+\Delta p'(y) ) =\int p(y)\ln\left(\frac{p(y)}{p(y)+\Delta p'(y)}\right)dy\;.\label{DF}
	\end{eqnarray}
	We then shall expand \eqref{DF} with respect to $p'$. Keeping only the first dominate term, we obtain
	\begin{eqnarray}
		\mathit{D}(p(y)\parallel p(y)+\Delta p'(y) ) \approx \int \frac{1}{2}\frac{\left(p'(y)\right)^2}{p(y)}dy=\frac{1}{2}I(p(y))\;.\label{DF2}
	\end{eqnarray}
	We could see that the right-hand side of \eqref{DF2} is nothing but the standard Fisher information. 
	\\
	\\
	Now we use two-parameter Kullback-Leibler divergence \cite{FracEn}
	\begin{equation}
		D_{q,q'}\left(p(y)\parallel p(y)+\Delta p'(y)\right)=\int p^q(y)\left(\ln\frac{p(y)}{p(y)+\Delta p'(y)}\right)^{q'}dy\;.
	\end{equation}
	Here we do again the expansion with respect to $p'$ and we obtain the two-parameter generalisation of the Fisher information from the first dominant term 
	\begin{equation}
		I_{a,b}(p)=\int p^a(y)\left(\frac{dp(y)}{dy}\right)^b dy\;, \label{2paraF}
	\end{equation}
	where $a=q-q'-1$ and $b=q'+1$ with the requirements $q > 0$ and $q'>0$. We find that, with a suitable choice of parameters, our whole hierarchy of Fisher information can be identified as shown in the table \ref{Ta2}.
	\begin{table}[htbp]
		\caption{Comparison our one-parameter Fisher information with two-parameter Fisher information \eqref{2paraF}. }
		\begin{center}
			\begin{tabular}{|c|c|c|c|}
				\hline
				\multirow{2}{5cm}{\centering $n^\text{th}$ Fisher information}& \multicolumn{2}{c|}{\centering Parameters} \\
				\cline{2-3} & \multicolumn{1}{p{1cm}|}{\centering $a$} & \multicolumn{1}{p{1cm}|}{\centering $b$}  \\ \hline
				$1^\text{st}$ order: $I_1$ & 1 & 2  \\
				$2^\text{nd}$ order: $I_2$ & 2 & 4  \\
				$3^\text{rd}$ order: $I_3$ & 3 & 6  \\
				$4^\text{th}$ order: $I_4$ & 4 & 8 \\
				\hline
			\end{tabular}
		\end{center} \label{Ta2}
	\end{table} 
	\\
	We note here that the quantities in \eqref{2paraF}, which are directly connected with our Fisher information hierarchy as we already mentioned, can be possibly viewed as the generalised Fisher information. However, there exist also other types of the generalised Fisher information for different purposes and motivations \cite{FANON2,GFM}.
	%
	%
	%
	%
	%
	%
	%
	
	\section{One-parameter generalised Fisher information matrix: one random variable}\label{section7}	
	With the result of the extension of the Fisher information in the case of the one estimated parameter, in this section, we will consider the case of the multi-estimated parameters: $\theta=\{\theta_1, \theta_2, ...,\theta_N\}$. Recalling the Fisher information matrix \eqref{Fismatrix}, we propose a one-parameter generalised Fisher information matrix
	\begin{eqnarray}
		[I_\lambda(\theta)] = \frac{1}{\lambda} \int dx \bigg(e^{\lambda  p(x \mid \theta ) [I(\theta)]}-\mathbb{I}_{N\times N}\bigg)\;,\label{FM1}
	\end{eqnarray}
	where $\mathbb{I}_{N\times N}$ is the identity matrix with dimension $N\times N$. Expanding \eqref{FM1} with respect to the parameter $\lambda$, one obtains
	\begin{eqnarray}
		[I_\lambda(\theta)] &=& \frac{1}{\lambda} \int dx \bigg( \lambda  p(x \mid \theta ) [I(\theta)] +
		\frac{\lambda^2}{2!}  p^2(x \mid \theta ) [I(\theta)]^2
		+ ...\bigg) \nonumber  \\ 
		&=& [I_1(\theta)]+ \frac{\lambda}{2!}[I_2(\theta)]+... \;. \label{Fmatrix}	
	\end{eqnarray}
	What we obtain now is a hierarchy of the Fisher information matrix: $\{[I_1(\theta)],[I_2(\theta)],...\}$. Obviously, $[I_1(\theta)]=[I(\theta)]$ is nothing but the standard Fisher information matrix. 
	\\
	\\
	It is not difficult to see that if the number of estimated parameters is one, \eqref{Fmatrix} simply reduces to
	\begin{eqnarray}
		I_\lambda(\theta) &=& \frac{1}{\lambda} \int dx \bigg(e^{[\partial_\theta \ln p]^2 }-1 \bigg) \nonumber \\
		&=&  \frac{1}{\lambda} \int dx \bigg( \lambda p [(\partial_\theta \ln p)^2] +
		\frac{\lambda^2}{2!} p^2 [(\partial_\theta \ln p)^2]^2
		+ ...\bigg) \;,\label{FM22}
	\end{eqnarray}
	where $p(x|\theta)=p$ and $\frac{\partial}{\partial \theta}=\partial_\theta $. Of course, \eqref{FM22} is what we obtained in the previous section and various properties are studied.
	\\
	\\
	Next, we would like to look further on the geometrical meaning of what we have in \eqref{FM1}. For simplicity, we shall focus on the case of two estimated parameters. Then we have
	\begin{eqnarray}
		[I_\lambda(\theta)] &=& \frac{1}{\lambda} \int dx \bigg(e^{\lambda p \begin{bmatrix}
				(\partial_1 \ln p)^2 & (\partial_1 \ln p)(\partial_2 \ln p) \nonumber \\
				(\partial_2 \ln p)(\partial_1 \ln p) & (\partial_2 \ln p)^2		
		\end{bmatrix}}
		-\mathbb{I}_{2\times 2}\bigg) \nonumber \\
		&=& \frac{1}{\lambda} \int dx \bigg( \lambda p 
		\begin{bmatrix}
			(\partial_1 \ln p)^2 & (\partial_1 \ln p)(\partial_2 \ln p) \\
			(\partial_2 \ln p)(\partial_1 \ln p) & (\partial_2 \ln p)^2		
		\end{bmatrix} \nonumber \\
		&+& 
		\frac{\lambda^2}{2!} p^2 {\begin{bmatrix}
				(\partial_1 \ln p)^2 & (\partial_1 \ln p)(\partial_2 \ln p) \\
				(\partial_2 \ln p)(\partial_1 \ln p) & (\partial_2 \ln p)^2		
		\end{bmatrix}}^2  
		+ ...\bigg) \nonumber \\
		&=& [I_1(\theta)]+ \frac{\lambda}{2!}[I_2(\theta)]+... \;, \label{FMM}
	\end{eqnarray}
	where $\{\partial_i=\partial/\partial\theta_i;i=1,2\}$.
	The first one is the standard Fisher information matrix with dimension $2\times2$ given by
	\begin{eqnarray}	
		[I_1(\theta)]&=&\begin{bmatrix}
			I_1(\theta)_{11}  & I_1(\theta)_{12}\\
			I_1(\theta)_{21} &I_1(\theta)_{22}		
		\end{bmatrix}
		=
		\begin{bmatrix}
			\int dx\;p(\partial_1 \ln p)^2 &  \int dx\;p(\partial_1 \ln p)(\partial_2 \ln p) \\
			\int dx\;p(\partial_2 \ln p)(\partial_1 \ln p) &  \int dx\; p(\partial_2 \ln p)^2		
		\end{bmatrix}\;. 	\label{FM21}
	\end{eqnarray}
	The 2$^\text{nd}$ order Fisher information matrix with the dimension $2\times2$ is given by
	\begin{eqnarray}	
		[I_2(\theta)]
		&=&  
		\begin{bmatrix}
			I_2(\theta)_{11}  & I_2(\theta)_{12}\\
			I_2(\theta)_{21} &I_2(\theta)_{22}		
		\end{bmatrix}\;,\label{FM2}	
	\end{eqnarray}
	where 
	\begin{subequations}
		\begin{alignat}{4}
			I_2(\theta)_{11} &=\int dx\; p^2 \bigg((\partial_i \ln p)^4+(\partial_i \ln p)^2(\partial_j \ln p)^2 \bigg)\;, \\            
			I_2(\theta)_{12}&=\int dx\; p^2 \bigg((\partial_i \ln p)^3(\partial_j \ln p)+(\partial_i \ln p)(\partial_j \ln p)^3\bigg)\;, \\ 
			I_2(\theta)_{21}&=\int dx\; p^2\bigg( (\partial_j \ln p)^3(\partial_i \ln p)+(\partial_j \ln p)(\partial_i \ln p)^3\bigg)\;,\\
			I_2(\theta)_{21}&=\int dx\; p^2 \bigg((\partial_i \ln p)^2(\partial_j \ln p)^2+(\partial_j \ln p)^4\bigg)\;.
		\end{alignat}
	\end{subequations}
	\\
	\textbf{Working example}: To illustrate the geometrical meaning of the first two Fisher information matrices in \eqref{FMM}, we shall again use the normal distribution 
	\begin{eqnarray}	
		p(x|\theta=\{\mu,\sigma\})=\frac{1}{\sqrt{2 \pi}\sigma} e^{-\frac{1}{2 \sigma^2}(x-\mu)^2} \;,
	\end{eqnarray}
	as a prime example. We therefore define a statistical manifold $M=\{N(\mu,\sigma^2) \}=\{p(x|\mu,\sigma)|(\mu,\sigma)\in \mathbb{R}\times \mathbb{R}^+\}$. It is not difficult to find that the log-likelihood is given by
	\begin{eqnarray}	
		\ln p(x|\theta=\{\mu,\sigma\})=\ln \frac{1}{\sqrt{2 \pi }\sigma} - \frac{1}{2 \sigma^2} (x-\mu)^2\;.
	\end{eqnarray}
	Therefore, we have
	\begin{eqnarray}	
		\frac{\partial}{\partial \theta_1}\ln p(x|\theta)=\frac{\partial}{\partial \mu}\ln p(x|\theta) = \frac{x-\mu}{\sigma^2}  \;\;\;
		\text{and}
		\;\;\;
		\frac{\partial}{\partial \theta_2}\ln p(x|\theta)=	\frac{\partial}{\partial \sigma}\ln p(x|\theta) =\frac{(x-\mu)^2}{\sigma^3} -\frac{1}{\sigma}\;.	
	\end{eqnarray}
	\\
	\emph{The first Fisher matrix $[I_1(\theta)]$}: We compute first the component $[I(\theta)]_{11}$ 
	\begin{eqnarray}	
		I_1(\theta)_{11} &=& \int_{-\infty}^{\infty} dx\;p(\partial_\mu \ln p)^2 \nonumber \\
		&=& \int dx \bigg( \frac{1}{\sqrt{2 \pi}\sigma} e^{-\frac{1}{2 \sigma^2}(x-\mu)^2} \bigg) \bigg( \frac{x-\mu}{\sigma^2} \bigg)^2 \nonumber \\
		&=& \frac{1}{\sqrt{2  \pi} \sigma^5} \int dz \; z^2 \bigg(e^{-\frac{1}{2 \sigma^2}z^2} \bigg)\;,\label{FM01}
	\end{eqnarray}
	where $z=x-\mu$ is used. Since \eqref{FM01} is the Gaussian integral, therefore, we obtain
	\begin{eqnarray}	
		I_1(\theta)_{11} &=& \frac{1}{ \sqrt{2\pi} \sigma^5}\biggl\{ \sqrt{2 \pi}  \sigma^3 \biggr\} 
		= \frac{1}{\sigma^2}\;.
	\end{eqnarray}
	We note here that this component of the Fisher information matrix is nothing but the standard one obtained in the previous section. Next, the component $I_1(\theta)_{22}$ can be computed as follows  
	\begin{eqnarray}	
		I_1(\theta)_{22} &=& \int dx\;p(\partial_\sigma \ln p)^2 \nonumber \\
		&=& \int dx \bigg( \frac{1}{\sqrt{2 \pi}\sigma} e^{-\frac{1}{2 \sigma^2}(x-\mu)^2} \bigg) \bigg( \frac{(x-\mu)^2}{\sigma^3} -\frac{1}{\sigma} \bigg)^2 \nonumber \\
		&=& \frac{1}{\sqrt{2  \pi} \sigma} \int dz  \bigg(e^{-\frac{1}{2 \sigma^2}z^2} \bigg) \bigg(\frac{z^4}{\sigma^6} -\frac{2z^2}{\sigma^4}+\frac{1}{\sigma^2} \bigg) \nonumber \\
		&=& \frac{1}{\sqrt{2  \pi} \sigma} 
		\bigg\{ \frac{1}{\sigma^6}\int dz\;z^4 \bigg(e^{-\frac{1}{2 \sigma^2}z^2} \bigg) - \frac{2}{\sigma^4} \int dz\;z^2 \bigg(e^{-\frac{1}{2 \sigma^2}z^2} \bigg) +  \frac{1}{\sigma^2} \int dz\; \bigg(e^{-\frac{1}{2 \sigma^2}z^2} \bigg)  \bigg\} \nonumber \\
		&=& \frac{1}{\sqrt{2  \pi} \sigma} 
		\bigg\{ \frac{1}{\sigma^6} \bigg( 3\sqrt{2\pi} \sigma^5 \bigg) - \frac{2}{\sigma^4} \bigg( \sqrt{2\pi} \sigma^3\bigg) +  \frac{1}{\sigma^2} \bigg( \sqrt{2 \pi} \sigma\bigg)  \bigg\} \nonumber \\
		&=& \frac{2}{\sigma^2} \;.
	\end{eqnarray}
	Lastly, we find that the off-diagonal components $I_1(\theta)_{12}$ and $I_1(\theta)_{21}$ are the same. These two quantities vanish as follows
	\begin{eqnarray}	
		I_1(\theta)_{12}=I_1(\theta)_{21} &=& \int dx\;p(\partial_\mu \ln p)(\partial_\sigma \ln p) \nonumber \\
		&=&  \int dx\; \bigg( \frac{1}{\sqrt{2 \pi}\sigma} e^{-\frac{1}{2 \sigma^2}(x-\mu)^2} \bigg) \bigg( \frac{x-\mu}{\sigma^2}\bigg) \bigg( \frac{(x-\mu)^2}{\sigma^3} -\frac{1}{\sigma} \bigg) \nonumber \\
		&=&  \frac{1}{\sqrt{2  \pi} \sigma} \int dz\; \bigg(e^{-\frac{1}{2 \sigma^2}z^2} \bigg) \bigg( \frac{z}{\sigma^2}\bigg) \bigg( \frac{z^2}{\sigma^3} -\frac{1}{\sigma} \bigg) \nonumber \\
		&=& \frac{1}{\sqrt{2 \pi} \sigma} \biggl\{  \frac{1}{\sigma^5}\int dz\;z^3 \bigg(e^{-\frac{1}{2 \sigma^2}z^2} \bigg) -  \frac{1}{\sigma^3}\int dz\;z \bigg(e^{-\frac{1}{2 \sigma^2}z^2} \bigg)  \biggr\} \nonumber \\
		&=& 0\;.
	\end{eqnarray}
	Finally, the standard Fisher information matrix can be expressed as 
	\begin{eqnarray}	
		[I(\theta)]=[I_1(\theta)]&=&
		\begin{bmatrix}
			\frac{1}{\sigma^2} & 0 \\
			0 &  \frac{2}{\sigma^2}		
		\end{bmatrix}\;. \label{SFmatrix}
	\end{eqnarray}
	Of course with \eqref{SFmatrix}, one can write the line element
	\begin{eqnarray}	
		ds^2= g_{ij} dx^{i} dx^{j}\;, \label{hmat}
	\end{eqnarray}
	where $i,j=1,2$. We now introduce a new set of variables such that 
	\begin{eqnarray}	
		\begin{bmatrix}
			dx^1 \\
			dx^2	
		\end{bmatrix}
		=\begin{bmatrix}
			d\mu \\
			d\sigma	
		\end{bmatrix}=\begin{bmatrix}
			dx \\
			dy	
		\end{bmatrix}\;. \label{gmat}
	\end{eqnarray}
	Then the matrix and its inverse become
	\begin{eqnarray}	
		g_{ij} = 
		\begin{bmatrix}
			1/y^2 & 0 \\
			0 & 1/y^2	
		\end{bmatrix}\;, \;\;\; 
		g^{ij} = 
		\begin{bmatrix}
			y^2 & 0 \\
			0 & y^2	
		\end{bmatrix}\;.
		\label{gmat2}
	\end{eqnarray}
	We note that the factor $2$ is dropped out at the moment for simplicity. With \eqref{gmat} and \eqref{gmat2}, \eqref{hmat} can be simply written as
	\begin{eqnarray}	
		ds^2= \frac{dx^2 + dy^2}{y^2}\;.
	\end{eqnarray}
	Now, the distance between two point $P$ and $Q$ in space $M$ is
	\begin{eqnarray}	
		s
		= \int_{P}^{Q} \mathscr{L}(y',y;x)dx\;,
	\end{eqnarray}
	where $y'=dy/dx$ and $\mathscr{L}(y',y;x)=\sqrt{1+y'^2/y^2}$ is treated as the Lagrangian. The Euler-Lagrange equation 
	\begin{eqnarray}	
		\frac{d}{dx} \frac{\partial \mathscr{L}(y',y;x)}{\partial y'} - \frac{\partial \mathscr{L}(y',y;x)}{\partial y}=0 \label{eleq}
	\end{eqnarray}
	yeilds 
	\begin{eqnarray}	
		yy''+(1+y'^2)=0\;. \label{eom1}
	\end{eqnarray}
	Equation \eqref{eom1} is the geodesic equation and the solution is given by
	\begin{eqnarray}	
		y&=&\sqrt{-2C_2 x- C_2^2+C_1-x^2} \nonumber \\
		&=& \sqrt{C_1 -(x+C_2)^2}\;, \label{g1}
	\end{eqnarray}
	where $C_1$ and $C_2$ are constants. Therefore, \eqref{g1} tells us that the shortest distance between $P$ and $Q$ will be on the semicircle-arc: $(x+C_2)^2+y^2=C_1$, where $(C_2,0)$ is the centre and $\sqrt{C_1}$ is the radius of the semicircle on the $M$. This manifold for the normal distribution is known as the hyperbolic half plane model or the Poincaré half plane model. In the case that we take back the factor $2$ into consideration, the solution will take the from
	\begin{eqnarray}	
		y&=&\frac{\sqrt{-2C_2 x- C_2^2+C_1-x^2}}{\sqrt{2}} \nonumber \\
		&=& \frac{\sqrt{C_1 -(x+C_2)^2}}{\sqrt{2}}\;,
	\end{eqnarray}
	which means that geodesics in $M$ is still half upper plane model, but the factor $\sqrt{2}$ curls it to be half-ellipses \cite{Suel}.
	\\
	\\
	Next, the curvature of the manifold $M$ will be computed. We then compute the first and second kind Christoffel symbols 
	\begin{eqnarray}	
		\Gamma_{\gamma \alpha \beta }=[\alpha \beta, \gamma]=\frac{1}{2} \bigg( \frac{\partial g_{\alpha \gamma} }{\partial x^\beta} + \frac{\partial g_{\beta \gamma} }{\partial x^\alpha} - \frac{\partial g_{\alpha \beta} }{\partial x^\gamma}\bigg)\;, \;\;\; \text{and} \;\;\; \Gamma^{\gamma}_{\alpha \beta} = g^{\gamma \delta}  	\Gamma_{\gamma \alpha \beta } \;,
	\end{eqnarray}
	resulting in
	\begin{eqnarray}	
		\Gamma^{1}_{12}=\Gamma^{1}_{21}=\Gamma^{2}_{22}= -1/\sigma \;, \;\;\; \text{and} \;\;\; \Gamma^{2}_{11}=1/2\sigma\;, 
	\end{eqnarray}
	while other components are zero. Many components of the Riemann curvature tensor 
	\begin{eqnarray}	
		R^{\eta}_{\alpha \beta \gamma} = \frac{\partial \Gamma^{\eta}_{\alpha \gamma}}{\partial x^\beta}-\frac{\partial \Gamma^{\eta}_{\beta \gamma}}{\partial x^\alpha}+\Gamma^{\tau}_{\alpha \gamma} \Gamma^{\eta}_{\tau \beta} - \Gamma^{\tau}_{\alpha \beta} \Gamma^{\eta}_{\tau \gamma}\;  
	\end{eqnarray}
	are zero and the remaining components are
	\begin{eqnarray}	
		R^{1}_{212}=-R^{1}_{221}=-1/\sigma^2\;\;\; \text{and} \;\;\;  	R^{2}_{121}=-R^{2}_{112}=-1/2\sigma^2\;.
	\end{eqnarray}
	The Ricci tensor can be computed by
	$	
	R^{\eta}_{\alpha \eta \gamma}=R_{\alpha \gamma}\;.
	$
	Again, all remaining components are 
	$	
	R_{11}=-1/2\sigma^2 $ and $ R_{22}=-1/\sigma^2\;.
	$
	Finally, we can find the Ricci scalar as
	\begin{eqnarray}	
		R=g^{\mu \nu} R_{\mu \nu}=\sigma^2(-1/2\sigma^2)+(\sigma^2/2)(-1/\sigma^2) = -1\;.
	\end{eqnarray}
	Then the curvature of the two-dimensional manifold $M$ is
	\begin{eqnarray}	
		\kappa = \frac{1}{2} R 
		= -\frac{1}{2}\;,
	\end{eqnarray}
	which implies that the manifold $M$ is the hyperbolic space, since its curvature is an uniform negative value $-1/2$. 
	\\
	\\
	\emph{The second Fisher matrix $[I_2(\theta)]$}: Next, we investigate the geometrical nature of the space $M$ with the matrix \eqref{FM2}. First, we compute
	
	\begin{eqnarray}	
		I_2(\theta)_{11}&=&\int dx\; p^2 \bigg((\partial_\mu \ln p)^4+(\partial_\mu \ln p)^2(\partial_\sigma \ln p)^2 \bigg) \nonumber \\
		&=& \int dx\; \bigg( \bigg( \frac{1}{2 \pi \sigma^2} e^{-\frac{1}{ \sigma^2}(x-\mu)^2} \bigg)  \bigg(\frac{x-\mu}{\sigma^2}\bigg)^4  \bigg) \nonumber\\
		&+& \int dx\; \bigg( \bigg( \frac{1}{2 \pi \sigma^2} e^{-\frac{1}{ \sigma^2}(x-\mu)^2} \bigg)  \bigg(\frac{x-\mu}{\sigma^2}\bigg)^2 \bigg(\frac{(x-\mu)^2}{\sigma^3}-\frac{1}{\sigma}\bigg)^2  \bigg) \nonumber \\
		&=&  \frac{1}{2 \pi \sigma^2} \bigg\{ \int dz\;  \bigg( e^{-\frac{1}{ \sigma^2}z^2} \bigg)  \bigg(\frac{z}{\sigma^2}\bigg)^4  
		+ \int dz\; \bigg(  e^{-\frac{1}{ \sigma^2}z^2} \bigg)  \bigg(\frac{z}{\sigma^2}\bigg)^2 \bigg(\frac{z^2}{\sigma^3}-\frac{1}{\sigma}\bigg)^2   \bigg\} \nonumber \\
		&=& \frac{1}{2 \pi \sigma^2} \biggl\{ \frac{1}{\sigma^{10}}\int dz\; z^6  \bigg(e^{-\frac{1}{ \sigma^2}z^2} \bigg) -\frac{1}{\sigma^8}\int dz\; z^4 \bigg(e^{-\frac{1}{ \sigma^2}z^2} \bigg)+\frac{1}{\sigma^6}\int dz\; z^2 \bigg(e^{-\frac{1}{ \sigma^2}z^2} \bigg)  \biggr\} \nonumber \\
		&=& \frac{1}{2 \pi \sigma^2} \biggl\{ \frac{1}{\sigma^{10}} \bigg(\frac{15 \sqrt{\pi} \sigma^7}{2^3} \bigg) -\frac{1}{\sigma^8}\bigg( \frac{3 \sqrt{\pi} \sigma^5}{2^2} \bigg)+\frac{1}{\sigma^6}\bigg( \frac{\sqrt{\pi}\sigma^3}{2} \bigg)  \biggr\} \nonumber \\
		&=& \frac{13\sqrt{\pi}}{2^4 \pi \sigma^5}\;, 
	\end{eqnarray}
	where we again use $z=x-\mu$. Similarly, the component $I_2(\theta)_{22}$ can be computed 
	\begin{eqnarray}	
		I_2(\theta)_{22}&=&\int dx\; p^2 \bigg((\partial_\sigma \ln p)^2(\partial_\mu\ln p)^2 +(\partial_\sigma \ln p)^4\bigg) \nonumber \\
		&=& \int dx\; \bigg( \bigg( \frac{1}{2 \pi \sigma^2} e^{-\frac{1}{ \sigma^2}(x-\mu)^2} \bigg)  \bigg(\frac{(x-\mu)^2}{\sigma^3}-\frac{1}{\sigma}\bigg)^2  \bigg(\frac{x-\mu}{\sigma^2}\bigg)^2  \bigg) \nonumber \\
		&+& \int dx\; \bigg( \bigg( \frac{1}{2 \pi \sigma^2} e^{-\frac{1}{ \sigma^2}(x-\mu)^2} \bigg) \bigg(\frac{(x-\mu)^2}{\sigma^3}-\frac{1}{\sigma}\bigg)^4  \bigg) \nonumber \\
		&=&  \frac{1}{2 \pi \sigma^2} \bigg\{ \int dz\;  \bigg( e^{-\frac{1}{ \sigma^2}z^2} \bigg)  \bigg(\frac{z^2}{\sigma^3} - \frac{1}{\sigma}\bigg)^2  \bigg(\frac{z}{\sigma^2}\bigg)^2  
		+ \int dz\; \bigg(  e^{-\frac{1}{ \sigma^2}z^2} \bigg)   \bigg(\frac{z^2}{\sigma^3}-\frac{1}{\sigma}\bigg)^4   \bigg\} \nonumber \\
		&=& \frac{1}{2 \pi \sigma^2} \biggl\{ \frac{1}{\sigma^{12}}\int dz\; z^8  \bigg(e^{-\frac{1}{ \sigma^2}z^2} \bigg) -\frac{3}{\sigma^{10}}\int dz\;z^6 \bigg(e^{-\frac{1}{ \sigma^2}z^2} \bigg) + \frac{4}{\sigma^{8}}\int dz\;z^4 \bigg(e^{-\frac{1}{ \sigma^2}z^2} \bigg) \nonumber \\
		&-& \frac{3}{\sigma^{6}}\int dz\;z^2 \bigg(e^{-\frac{1}{ \sigma^2}z^2} \bigg)  +\frac{1}{\sigma^{4}}\int dz\; \bigg(e^{-\frac{1}{ \sigma^2}z^2} \bigg) \biggr\} \nonumber \\
		&=& \frac{1}{2 \pi \sigma^2} \biggl\{ \frac{1}{\sigma^{12}} \bigg( \frac{7 \cdot 5 \cdot 3 \sqrt{\pi} \sigma^9 }{2^4}\bigg) -\frac{3}{\sigma^{10}} \bigg(  \frac{ 5 \cdot 3 \sqrt{\pi} \sigma^7 }{2^3} \bigg) + \frac{4}{\sigma^{8}} \bigg( \frac{3 \sqrt{\pi} \sigma^5 }{2^2} \bigg) \nonumber \\
		&-& \frac{3}{\sigma^{6}} \bigg(  \frac{\sqrt{\pi} \sigma^3 }{2} \bigg)  +\frac{1}{\sigma^{4}} \bigg(\sqrt{\pi} \sigma\bigg) \biggr\} \nonumber \\
		&=& \frac{55\sqrt{\pi}}{2^5 \pi \sigma^5}
		\;. 
	\end{eqnarray}
	For the off-diagonal components, $I_2(\theta)_{12}$ and $I_2(\theta)_{21}$, they are all zero as follows
	\begin{eqnarray}	
		I_2(\theta)_{12}=I_2(\theta)_{21}&=&\int dx\; p^2 \bigg((\partial_\mu \ln p)^3(\partial_\sigma \ln p)+(\partial_\mu \ln p)(\partial_\sigma \ln p)^3\bigg) \nonumber \\
		&=& \int dx\; \bigg( \bigg( \frac{1}{2 \pi \sigma^2} e^{-\frac{1}{ \sigma^2}(x-\mu)^2} \bigg) \bigg(\frac{x-\mu}{\sigma^2}\bigg)^3 \bigg(\frac{(x-\mu)^2}{\sigma^3}-\frac{1}{\sigma}\bigg) \bigg) \nonumber \\
		&+& \int dx\; \bigg( \bigg( \frac{1}{2 \pi \sigma^2} e^{-\frac{1}{ \sigma^2}(x-\mu)^2} \bigg) \bigg(\frac{x-\mu}{\sigma^2}\bigg)  \bigg(\frac{(x-\mu)^2}{\sigma^3}-\frac{1}{\sigma}\bigg)^3  \bigg) \nonumber \\
		&=& \frac{1}{2 \pi \sigma^2} \bigg\{ \int dx\;  \bigg(  e^{-\frac{1}{ \sigma^2}z^2} \bigg) \bigg(\frac{z}{\sigma^2}\bigg)  \bigg(\frac{z^2}{\sigma^3}-\frac{1}{\sigma}\bigg)^3 \nonumber \\
		&+&  \int dx\;  \bigg(  e^{-\frac{1}{ \sigma^2}z^2} \bigg) \bigg(\frac{z}{\sigma^2}\bigg)^3  \bigg(\frac{z^2}{\sigma^3}-\frac{1}{\sigma}\bigg) \bigg\} \nonumber \\
		&=& \frac{1}{2 \pi \sigma^2} \biggl\{ \frac{1}{\sigma^{11}}\int dz\; z^7 \bigg(e^{-\frac{1}{ \sigma^2}z^2} \bigg) +\frac{2}{\sigma^9}\int dz\;z^5 \bigg(e^{-\frac{1}{ \sigma^2}z^2} \bigg) \nonumber \\
		&+& \frac{2}{\sigma^7}\int dz\;z^3 \bigg(e^{-\frac{1}{ \sigma^2}z^2} \bigg) +\frac{1}{\sigma^5}\int dz\;z \bigg(e^{-\frac{1}{ \sigma^2}z^2} \bigg) \biggr\} \nonumber \\
		&=& 0\;. 
	\end{eqnarray}
	Therefore, for the normal distribution, the matrix \eqref{FM2} can be explicitly expressed as
	\begin{eqnarray}	
		[I_2(\theta)]&=&
		\begin{bmatrix}
			\frac{13\sqrt{\pi}}{2^4 \pi \sigma^5} & 0 \\
			0 & \frac{55\sqrt{\pi}}{2^5 \pi \sigma^5}	
		\end{bmatrix}\equiv \Tilde{g}_{ij}\;. \label{2ndFM}
	\end{eqnarray}
	The element length in this case becomes
	\begin{eqnarray}	
		ds^2= \Tilde{g}_{ij} dx^{i} dx^{j}   \;\;\; \text{or} \;\;\;  ds^2= \frac{Adx^2 + Bdy^2}{y^5}\;,\;\text{where}\;\;
		\begin{bmatrix}
			dx^1 \\
			dx^2	
		\end{bmatrix}
		=\begin{bmatrix}
			d\mu \\
			d\sigma	
		\end{bmatrix}=\begin{bmatrix}
			dx \\
			dy	
		\end{bmatrix}\;. \label{Fmat}
	\end{eqnarray}
	Here $A=\frac{13\sqrt{\pi}}{2^4 \pi }$ and $B=\frac{55\sqrt{\pi}}{2^5 \pi }$. 
	The distance between two point $P$ and $Q$ in space $M$ can be defined by 
	\begin{eqnarray}	
		s
		&=& \int_{P}^{Q} \mathscr{L}(y',y;x)dx \;,
	\end{eqnarray}
	where $\mathscr{L}(y',y;x)=\sqrt{\frac{A+By'^2}{y^5}}$ is treated as the Lagrangian. We notice that, to obtain the geodesic equation, it is not convenient to use the Euler-Lagrange equation in the form \eqref{eleq}. With this particular form of the Lagrangian, it is better to use the alternative form for Euler-Lagrange equation 
	\begin{eqnarray}	
		\frac{d}{dx} \bigg( \mathscr{L}(y',y;x)-y' \frac{\partial \mathscr{L}(y',y;x)}{\partial y'}  \bigg) -\frac{\partial \mathscr{L}(y',y;x)}{\partial x}=0\;, \label{elqal}
	\end{eqnarray}
	which is applicable if $\mathscr{L}$ is not an explicit function of variable $x$.
	Hence, we obtain
	\begin{eqnarray}	
		y'=\sqrt{\frac{D_2}{y^5}-D_3}\;,
	\end{eqnarray}
	where $D_2=\frac{A^2}{B D_1^2}$, $D_3=\frac{A}{B}$ and $D_1$ is a constant of integration. The solution, in terms of the original variables, is given by 
	\begin{eqnarray}	
		x= 2 y^{7/2} \left(\frac{\sqrt{1-\frac{D_3 y^5}{D_2}}}{7 \sqrt{D_2-D_3 y^5}}\right) {}_2 F_1 \left(\frac{1}{2}, \frac{7}{10}; \frac{17}{10}; \frac{D_3 y^5}{D_2}\right)-E\;,
	\end{eqnarray}
	where ${}_2 F_1 \left(a, b; c; z \right)$ is the hypergeometric function and $E$ is another constant of integration. 
	The non-zero components of the Christoffel symbols in this case are given by
	\begin{eqnarray}	
		\Gamma^{1}_{12}=\Gamma^{1}_{21}=\Gamma^{2}_{22}= -5/2 y \;, \;\;\; \text{and} \;\;\; \Gamma^{2}_{11}=5A/2B y\;,
	\end{eqnarray}
	the non-zero components of the Riemann curvature tensor are given by
	\begin{eqnarray}	
		R^{1}_{212}=-R^{1}_{221}=-5/2 y^2\;\;\; \text{and} \;\;\;  	R^{2}_{121}=-R^{2}_{112}=-5A/2 y^2\;.
	\end{eqnarray}
	The Ricci tensor survives only two components
	\begin{eqnarray}	
		R_{11}=-5A/2B y^2 \;\;\; \text{and} \;\;\; R_{22}=-5/2 y^2\;.
	\end{eqnarray}
	Finally, we can find that the Ricci scalar is
	\begin{eqnarray}	
		R=\Tilde{g}^{ij} R_{ij}=(y^5/A)(-5/2B y^2)+(y^5/B)(-5/2 y^2) = -\bigg(5(1+A)/2AB\bigg)y^3\;,
	\end{eqnarray}
	and then the curvature becomes
	\begin{eqnarray}	
		\kappa = \frac{1}{2} R 
		= \frac{1}{2} \left(-\frac{5(1+A)}{2AB} y^3\right)
		= -\frac{5(1+A)}{4AB} y^3\;.
	\end{eqnarray}
	Here, what we can conclude is that the curvature of the space $M$, equipped with the matrix \eqref{FM2}, is not uniform since it depends solely on the variable $y$ (which is actually the variable $\sigma$).
	\section{Conclusion}\label{section8}
	We succeed to construct the one-parameter generalised Fisher information in the case of the one random variable and one estimated parameter. The main method used to derive the one-parameter generalised Fisher information is the variational principle. We consider here the Fisher information as the action functional of the free particle Lagrangian. With the new insight of the one-parameter generalised Lagrangian \cite{MulLag}, one can naturally obtain the one-parameter generalised Fisher information. Furthermore, we can treat our one-parameter generalised Fisher information as the generator to obtain what we call the Fisher information hierarchy. The generalised Cramér-Rao inequality is also obtained with the help of the Hölder's inequality. We find that our Fisher information hierarchy does not follow the addition property, expect for the standard Fisher information. The interesting point is that this non-additive property pops up when the Tsallis index is not equal to 1 as shown in table \ref{Ta1}. Of course, this point is quite interesting since this non-additive property has been widely discussed in the Tsallis statistic \cite{Tsallis}. 
	Let us point out possibly indirect connection with the Tsallis entropy by recalling our one-parameter generalised Fisher information 
	\[
	I_\lambda[p]=\frac{4}{\lambda}\int \left[e^{\frac{\lambda}{4} \frac{p'^2(y)}{p(y)}}-1 \right]dy\;,
	\]
	and Tsallis entropy 
	\[
	S_q[p]=\frac{1}{1-q}\left[\int p^q(y)-1 \right]dy\;.
	\]
	It might be seem a bit strange but these two quantities more or less similar in the sense that they both contain a parameter and under the suitable limit the standard quantities can be recovered. However, more direction connection with the entropy might be the relative entropy or the two-parameter Kullback–Leibler divergence and our whole hierarchy Fisher information can be identified with the two-parameter Fisher information with the appropriate choice of parameters. 
	\\
	\\
	In the case of the one random variable with multi-estimated parameters, the one-parameter generalised Fisher information matrix is obtained. The hierarchy of Fisher information matrices is computed. Geometrical interpretation for the first two Fisher information matrices in the hierarchy is investigated with well known example: the normal distribution. With the normal distribution, the statistical manifold is two-dimensional upper plane. The standard Fisher matrix gives a constant negative curvature and therefore, the space is hyperbolic. The geodesic between two points on the space will be on the ellipse. The $2^{\text{nd}}$ order Fisher information matrix gives different geometrical nature on this very same space. The curvature is not constant and the geodesic equation is in the form of the hypergeometric function. 
	\\
	\\
	We finally note that the application of the one-parameter generalised Fisher information and one-parameter generalised Fisher information matrix is not so obvious and we shall leave this as an open problem for further investigation.
	\section*{Acknowledgments}
	Worachet Bukaew would like to acknowledge the partial financial support provided by The Prof. Dr. Sujin Jinahyon Foundation. Authors are very indebted to Dr. Ninnat Dangniam for giving us an interested point on the Fisher information matrix part. Additionally, authors are grateful to Mr. Amornthep Tita for his suggestion on computation of the curvature.
	
\end{document}